\documentclass[a4paper,11pt]{article}

\usepackage{enumitem}
\usepackage[margin=1in]{geometry}

\usepackage{natbib}
 \bibpunct[, ]{(}{)}{,}{a}{}{,}%

\usepackage{makecell}
\usepackage{fancyhdr}
\usepackage{amsmath}
\usepackage{amsthm}
\usepackage{amsthm,amssymb}
\usepackage{graphicx}
\usepackage{caption}
\usepackage{authblk}

\newtheorem{definition}{Definition}
\newtheorem{proposition}{Proposition}

\title{New Exact Algorithm and Solution Properties for the Vehicle Routing Problem with Stochastic Demands}
\author[1]{Alexandre Florio\thanks{alexandre.de.macedo.florio@univie.ac.at}}
\author[1]{Richard Hartl\thanks{richard.hartl@univie.ac.at}}
\author[2]{Stefan Minner\thanks{stefan.minner@tum.de}}
\affil[1]{Department of Business Administration, University of Vienna}
\affil[2]{Logistics and Supply Chain Management, Technical University of Munich}

\begin{document}

\maketitle

\abstract{This paper considers the vehicle routing problem with stochastic demands (VRPSD) under optimal restocking. We develop an exact algorithm that is effective for solving instances with many vehicles and few customers per route. In our experiments, we show that in these instances solving the stochastic problem is most relevant (i.e., the potential gains over the deterministic equivalent solution are highest). The proposed branch-price-and-cut algorithm relies on an efficient labeling procedure, exact and heuristic dominance rules, and completion bounds to price profitable columns. Instances with up to 76 nodes could be solved in less than 5 hours, and instances with up to 148 nodes could be solved in long-runs of the algorithm. The experiments also allowed new findings on the problem. Solving the stochastic problem leads to solutions up to 10\% superior to the deterministic equivalent solution. When the number of routes is not fixed, the optimal solutions under detour-to-depot and optimal restocking are nearly equivalent. Opening new routes is a good strategy to reduce restocking costs, and in many cases results in solutions with less transportation costs. For the first time, scenarios where the expected demand in a route is allowed to exceed the capacity of the vehicle were also tested, and the results indicate that superior solutions with lower cost and fewer routes exist.}

\section{Introduction}
% problem introduction
We consider a vehicle routing problem \citep[VRP, see][]{vrppma} in which the demands of the customers are stochastic, and only disclosed upon arrival of the vehicle at their locations. This problem is known as the vehicle routing problem with stochastic demands (VRPSD), and it is the most studied stochastic variant of the VRP \citep{Gendreauetal2016}. In the VRPSD, each vehicle is allowed to perform replenishment trips to the depot, in order to increase the total capacity available for a route. These restocking trips might be \emph{reactive}, i.e., triggered by a stockout when serving the demand of some customer, or \emph{preventive}, i.e., performed in anticipation of a probable stockout later in the route. Therefore, after every visit a decision has to be made: whether or not to replenish before visiting the next customer in the route. It is well-known that, given an a priori route, optimal restocking decisions can be computed with a stochastic dynamic programming algorithm \citep{YeeGolden1980}. We consider the VRPSD under optimal restocking, which calls for the identification of a set of a priori routes visiting all customers with minimum total expected cost, considering that restocking decisions are always optimal.

% HERE
Over the years, many exact algorithms have been developed for the VRPSD under the reactive, or \emph{detour-to-depot}, policy. Under this policy, the vehicle cannot perform preventive restocking trips, even if there is a high stockout probability at some customer located far from the depot. Recently, the first approaches addressing the VRPSD under optimal restocking have been proposed \citep{louveaux2017,salavati2017exact}. These branch-and-cut algorithms are effective for solving instances with many customers per route (from 10 up to 50). We show that in many of these instances, the optimal stochastic solution is only marginally better than the deterministic equivalent solution. Traditional (deterministic) VRP instances are often characterized by much smaller customers per route ratios, and are likely to be intractable to the current methods. However, as shown in our computational experiments, these are the instances where solving the stochastic problem is more meaningful, i.e., where the value of the stochastic solution tends to be higher. Therefore, an exact algorithm for the VRPSD under optimal restocking, which can tackle the instances where solving the stochastic problem is mostly relevant, is still an open question in the literature.

% our methodological approach
We close this gap by designing a branch-price-and-cut algorithm for the VRPSD. The main component of our algorithm is an efficient labeling procedure for pricing profitable columns. This procedure leverages on stochastic dynamic programming for route cost evaluation. Therefore, no assumptions other than independency and discreteness are placed on the probability distributions of the demands. In fact, the demands of different customers could follow different probability distributions. Exact and approximate dominance rules are proposed and applied to control the combinatorial growth of labels, together with reduced cost (or completion) bounds. We also use a few techniques successfully employed by state-of-the-art algorithms for the deterministic capacitated VRP: \emph{ng}-route relaxation \citep{baldacci2011new}, and separation of rounded capacity cuts.

% statement of contributions (load analysis, VSS, method, literature instances, dominance rules)
In summary, in this paper we report the following contributions and findings:

\begin{enumerate}[label={(\arabic*)}]
\item An algorithm for the VRPSD under optimal restocking and general demand distributions. This is the first algorithm for the problem that is effective for solving instances with many vehicles and few customers per route. As we show, in these instances solving the stochastic problem is particularly relevant. Instances of moderate size (up to 75 nodes) can be solved in less than 5 hours, and larger instances (up to 148 nodes) can be solved in long-runs of the algorithm;
\item Computation of the value of the stochastic solution of all instances solved, which shows that the optimal stochastic solution can be up to 10.5\% superior to the deterministic equivalent solution. The potential gains are higher on instances with fewer customers per route, which are also the instances where the proposed algorithm performs better. On many instances with an average of 15 or more customers per route, solving the stochastic problem offers little to no value, at least under the probability distributions considered;
\item A relaxed label dominance rule, which accelerates the branch-price-and-cut algorithm, and also enables it to run in heuristic mode. This heuristic performs remarkably well, finding the optimum in every instance where the optimum is known, usually in much less computational time;
\item Comparison of the detour-to-depot and the optimal restocking policies, showing that the optimal policy offers little to no advantage over the detour-to-depot policy, when the number of routes is not fixed and the total expected demand in a route cannot exceed the capacity of the vehicle, i.e., when the load factor is up to 1;
\item Finally, experiments with load factors larger than 1 demonstrate that the overall costs can be reduced, when the expected demand in a route is allowed to exceed the capacity of the vehicle, even if only slightly. This is the first time different load scenarios have been considered in the (multi-vehicle) VRPSD.
\end{enumerate}

The remainder of this paper is organized as follows: in \S\ref{literature}, we review the main contributions on exact methods for the VRPSD; in \S\ref{formulations}, we present the VRPSD formulations that serve as the basis for the branch-price-and-cut algorithm, which is described in \S\ref{algorithm}. The results are presented and discussed in \S\ref{results}, and conclusions and directions for future research in \S\ref{conclusions}.

\section{Related Literature} \label{literature}
The first work on the VRPSD is attributed to \cite{tillman1969multiple}, who proposed a simple heuristic for solving the problem in a multi-depot setting. Since then, the problem has been receiving increasing attention from researchers, who over the years developed exact and heuristic approaches for solving it. We focus our review on exact methods for the VRPSD. A review on approximate approaches, and also on other stochastic variants of the VRP, can be found in \cite{gendreau2014stochastic}.

% recourse actions, restocking policies, detour-to-depot and optimal policy
When compared to the deterministic VRP, a key difference in the VRPSD is that the total demand in a route may exceed the capacity of the vehicle. There are a number of alternatives for dealing with this possibility. Early methods considered chance-constrained or penalty-based models \citep{Droretal1993,Laporteetal1989}. Nevertheless, allowing restocking trips is the alternative that has been given more attention by researchers. In fact, by permitting such trips we can still enforce that all customers have their entire demand fulfilled, as in the VRP.

% HERE
When replenishment is allowed, one must also specify the restocking policy adopted. The simple detour-to-depot policy, first stated in \cite{Droretal1989}, allows replenishment only when the remaining quantity in the vehicle is not sufficient to serve the current customer, an event known as stockout, or \emph{route failure}. On the other hand, the optimal policy prescribes replenishments in an optimal way, and can be computed with the stochastic dynamic programming algorithm from \cite{YeeGolden1980}. In optimal restocking, preventive replenishment can be performed to avoid route failures on customers located far from the depot. The advantages of optimal restocking, however, only materialize in routes where the probability of the actual demand exceeding the capacity of the vehicle is significant, i.e., in routes where replenishment may actually be needed. A comparison of both policies can be found in \cite{florio2017}.

Under the detour-to-depot policy, the VRPSD was solved for the first time in \cite{Gendreauetal1995}. The algorithm was based on the integer L-shaped method \citep{LaporteLouveaux1993}, a general branch-and-cut procedure for solving two-stage stochastic programs with binary first-stage variables. The method was also used in the exact approaches from \cite{Laporteetal2002} and \cite{Jabalietal2014}. The first approach based on branch-and-price was proposed by \cite{ChristiansenLysgaard2007}. It was significantly improved by \cite{Gauvinetal2014}, who added to the algorithm separation of valid inequalities (transforming it into a branch-price-and-cut) and implemented techniques already used successfully for solving the deterministic VRP.

% HERE
Under the optimal policy, the first approach in the literature is the one from \cite{louveaux2017}. The proposed algorithm, based on the integer L-shaped method, could solve instances with up to 101 nodes. Two methods for bounding the expected restocking cost were developed, and computational experiments were conducted assuming identically distributed customer demands. In \cite{salavati2017exact} new bounds were introduced, and instances with up to 60 customers and non-identically distributed demands could be solved. When replicating some of the results in \cite{louveaux2017}, we found out that, in many of the instances, the optimal stochastic solution is only marginally (less than 1\%) better than the deterministic equivalent solution. Thus, identifying the characteristics of VRPSD instances where solving the stochastic problem is more relevant is also an open question.

Another assumption in both these approaches is the consideration of a maximum load factor of 1, i.e., constraints limiting the maximum expected demand in a route to be not greater than the capacity of the vehicle. While such constraint may be reasonable in a detour-to-depot setting, where a replenishment can be potentially costly, it restricts unnecessarily the feasible solution space when optimal restocking is applied. As we show, better solutions can be found when this constraint is relaxed, even if only slightly, which is in line with the known result that, in the VRPSD under optimal restocking, there is always one optimal solution where a single vehicle visits all the customers \citep{Yangetal2000}.

Regarding the instances solved by the above methods, we note that only the branch(-price)-and-cut approaches from \cite{ChristiansenLysgaard2007} and \cite{Gauvinetal2014} are capable of solving literature instances (under the detour-to-depot policy) with no adaptations to the vehicle capacity. The experiments performed in \cite{louveaux2017} consider depot and customer locations according to 7 literature instances, but the vehicle capacity is modified to enable optimal solutions with up to only 3 vehicles. The other algorithms based on the integer L-shaped method perform most of the experiments on randomly generated instances.

When a single vehicle is considered, a mixed-integer linear model for the problem can be derived using Markov decision theory. The model presented in \cite{florio2017} computes, simultaneously, optimal a priori tours and restocking policies, and is valid for general discrete probability distributions and load factors. Unfortunately, the dimension of the model grows polynomially with the capacity of the vehicle (so, in fact, exponentially with the size of the instance). For this reason, only simplified instances can be solved, which is enough for measuring the sub-optimality of the detour-to-depot policy, but not enough for solving practical instances of the VRPSD.

\section{VRPSD Formulations} \label{formulations}
The branch-price-and-cut algorithm presented in \S\ref{algorithm} separates cuts defined over the variables of a compact formulation, and computes the linear bound by applying column generation on an extended formulation. In this section, we present both formulations.

In what follows, the VRPSD is defined on a graph with $N+1$ nodes. Node $v_{0}$ denotes the depot, and nodes $v_{1},\ldots,v_{N}$ the customers to be served. The demand of customer $v_{i}$ is modeled by the discrete random variable $D_{i}$, with probability mass function $p_{D_{i}}(d)$, which we assume has a (non-negative) integer support. We also assume the demands of the customers are independent. Let $d_{ij}$ be the transportation cost when traveling from node $v_{i}$ to node $v_{j}$. The triangle inequality holds. A number of vehicles is available at the depot, each with a positive and integer capacity of $Q$. A feasible route starts at the depot, visits one or more customers, and returns to the depot. The \emph{route load} is the sum of the expected demands of the customers visited in the route. The route load cannot exceed $fQ$, where $f$ is the positive and real-valued load factor parameter.

\subsection{Compact Formulation}
Using binary variables $x_{ij}$ to indicate whether node $v_{j}$ is visited immediately after node $v_{i}$, the VRPSD can be formulated as the following two-stage stochastic program:
\begin{align}
&\text{minimize}&&\sum_{i=0}^{N}\sum_{i=0}^{N}d_{ij}x_{ij}+E\left[\mathcal{Q}(\mathbf{x})\right]\,, \nonumber \\
&\text{subject to}&&\sum_{i=0}^{N}x_{ij}=1&\qquad\forall j\in\{1,\ldots,N\}\,, \nonumber \\
&&&\sum_{j=0}^{N}x_{ij}=1&\qquad\forall i\in\{1,\ldots,N\}\,, \nonumber \\
&&&\sum_{v_{i}\in S}\sum_{v_{j}\notin S}x_{ij}\geq\left\lceil\sum_{v_{i}\in S}E[D_{i}]/fQ\right\rceil&\qquad\forall S\subset\{v_{1},\ldots,v_{N}\}\,, \label{rcc} \\
&&&x_{ij}\in\{0,1\}\,, \nonumber
\end{align}
where $\mathbf{x}=\{x_{ij}:i,j\in\{0,\ldots,N\}\}$ is the set of first-stage decision variables, $\mathcal{Q}(\mathbf{x})$ is a random variable representing the \emph{restocking cost} when the first-stage decision is $\mathbf{x}$ and the optimal restocking policy is applied, and $E[.]$ is the expectation operator. We call the first term in the objective function the \emph{a priori cost}, which is deterministic once a first-stage decision has been made, i.e., once the a priori routes have been defined. We are not interested in further defining $\mathcal{Q}(\mathbf{x})$, since the evaluation of the objective value is done according to the formulation and procedures described in \S\ref{extended}.

Constraints \eqref{rcc} are the so-called \emph{rounded capacity cuts}, and are already adapted to allow routes with loads up to $fQ$. The separation of these cuts is a \emph{NP}-hard problem, but efficient heuristics are available. In many algorithms for vehicle routing problems, these cuts are separated using the heuristic procedures implemented in the $\mathtt{CVRPSEP}$ package \citep{lysgaard2003cvrpsep}. We do the same in our implementation.

\subsection{Extended Formulation} \label{extended}
A feasible route can be represented as a sequence $r=\{v_{0}, v_{s_{1}}, \ldots, v_{s_{l}}, v_{0}\}$, such that $l\geq 1$, $s_{i}\in\{1,\ldots,N\}$, $s_{i}\neq s_{j}$ if $i\neq j$, and:
\begin{equation} \nonumber
\sum_{v_{i}\,:\,r, v_{i}\neq v_{0}}E[D_{i}]\leq fQ\,,
\end{equation}
where the summation is over every element of $r$ except $v_{0}$.

The demand of a customer is only disclosed upon arrival of the vehicle at the customer location. If the remaining quantity in the vehicle is not sufficient to serve the entire demand, then the available quantity is delivered to the customer, and the vehicle performs replenishment trips to the depot until the entire demand of the customer is satisfied. After serving some customer, the vehicle is allowed to perform a replenishment trip to the depot before visiting the next customer in the a priori sequence. We assume the restocking decisions are made optimally. Therefore, the expected cost of a route can be calculated with the stochastic dynamic programming algorithm from \cite{YeeGolden1980}. The version of the algorithm we present below allows unbounded demands, i.e., for some customer $v_{i}$ and some $d>Q$, it is possible that $p_{D_{i}}(d)>0$. This is particularly important, since we experiment with unbounded distributions (e.g., Poisson), in which the probability of demands above $Q$ may be significant, at least for some customers. Let $r=\{v_{0}, v_{s_{1}}, \ldots, v_{s_{l}}, v_{0}\}$ be a feasible route. We define $s_{0}=0$. The expected cost of route $r$, denoted by $c_{r}$, can be calculated as follows:
\begin{equation} \nonumber
c_{r}=\nu^{r}(0,Q)\,,
\end{equation}
where $\nu^{r}(i,q)$ is the expected remaining cost of route $r$, given that the vehicle has just served the entire demand of customer $v_{s_{i}}$ (or is at the depot ready to start the route if $i=0$), and has a remaining quantity of $q$:
\begin{equation} \nonumber %\label{eq:defviq}
    \nu^{r}(i,q)=
\begin{cases}
    \min\{\nu_{dir}^{r}(i,q), \nu_{rep}^{r}(i)\},& \text{if } i<l \,, \\
    d_{s_{l}0},& \text{if } i=l \,.
\end{cases}
\end{equation}

The values $\nu_{dir}^{r}(i,q)$ and $\nu_{rep}^{r}(i)$ correspond, respectively, to the expected costs of the decisions of going to the next customer directly, or after a replenishment trip to the depot:
\begin{align}
\nu_{dir}^{r}(i,q)&=d_{s_{i}s_{i+1}}+\sum_{k=0}^{\infty}\left[(d_{s_{i+1}0}+d_{0s_{i+1}})\Psi(k,q)+\nu^{r}(i+1, Q\Psi(k,q)+q-k)\right]p_{D_{s_{i+1}}}(k) \,, \nonumber \\
\nu_{rep}^{r}(i)&=d_{s_{i}0}+d_{0s_{i+1}}+\sum_{k=0}^{\infty}\left[(d_{s_{i+1}0}+d_{0s_{i+1}})\Psi(k,Q)+\nu^{r}(i+1, Q\Psi(k,Q)+Q-k)\right]p_{D_{s_{i+1}}}(k) \,. \nonumber
\end{align}

In these expressions, the term $\Psi(k,q)$ gives the number of replenishment trips necessary to fulfill a demand of $k$, given that the remaining quantity in the vehicle is $q$:
\begin{equation} \nonumber
    \Psi(k, q)=
\begin{cases}
    0,& \text{if } k\leq q\,, \\
    1+\Psi(k, q+Q),& \text{otherwise} \,.
\end{cases}
\end{equation}

Or, in closed-form expression:
\begin{equation} \nonumber
\Psi(k, q)=\left\lceil\frac{k-q}{Q}\right\rceil^+ \,.
\end{equation}

Let $\mathcal{R}$ be the set of all feasible routes to some VRPSD instance. In addition, for some route $r\in\mathcal{R}$, let $a_{ir}$ be a binary value indicating whether customer $v_{i}\in r$, i.e., whether route $r$ visits customer $v_{i}$. Then, the VRPSD can be formulated as the following set-partitioning problem:
\begin{align}
&\text{minimize}&&\sum_{r\in\mathcal{R}}c_{r}\theta_{r}\,, \nonumber \\
&\text{subject to}&&\sum_{r\in\mathcal{R}}a_{ir}\theta_{r}=1&\qquad\forall i\in\{1,\ldots,N\}\,, \label{spc} \\
&&&\theta_{r}\in\{0,1\}\,. \nonumber
\end{align}

If the number of routes is fixed to some value $U$, an additional constraint $\sum_{r\in\mathcal{R}}\theta_{r}=U$ must be observed. When rounded capacity cuts are added to the extended formulation, they are translated into the following additional inequalities \citep[see][]{desaulniers2011cutting}:
\begin{align}
&&&\sum_{r\in\mathcal{R}}\sum_{v_{i}\in S_{k}}\sum_{v_{j}\notin S_{k}}b_{r}^{i,j}\theta_{r}\geq\left\lceil\sum_{v_{i}\in S_{k}}E[D_{i}]/fQ\right\rceil&\forall k\in\{1,\ldots,C\} \label{rccsp} \,,
\end{align}
where $C$ is the number of cuts added, $S_{k}\subset\{v_{1},\ldots,v_{N}\}$, $k\in\{1,\ldots,C\}$, are the sets that define each cut, and $b_{r}^{i,j}$ is a binary value indicating whether node $v_{j}$ is visited immediately after node $v_{i}$ in $r$.

Our decision of including only rounded capacity cuts has two main reasons. First, in the extended formulation of the capacitated VRP, these are the inequalities that yield the highest improvement in the linear bound, among the families of cuts that can be defined in terms of the variables of the compact formulation \citep{fukasawa2006robust}. We expect this result to also hold for the VRPSD, since in both problems the feasible regions are related. Second, adding other classes of cuts \citep[e.g., non-robust cuts, see][]{poggi2014new} would complicate even more the already involved pricing problem that we describe next.

\section{Branch-Price-and-Cut Algorithm} \label{algorithm}
In many variants of the VRP, set-partitioning formulations give relatively tight linear bounds, and are thus appealing for use within a branch-and-bound framework. Nevertheless, the exponential number of variables means that the linear bound must be computed by delayed column generation \citep[see][]{Desaulniersetal2005,lubbecke2005selected}. The resulting technique is called branch-and-price \citep{barnhart1998branch}, or branch-price-and-cut when valid inequalities are added during the execution of the algorithm.

In branch-and-price approaches for vehicle routing problems, the main difficulty is usually the solution of the pricing problem, i.e., the problem of finding a profitable variable, when solving the linear relaxation of the extended formulation by column generation. In our case this difficulty is magnified, since the cost of a route is computed by stochastic dynamic programming. The main components of our algorithm are the labeling procedures, dominance rules and completion bounds that enable the efficient solution of the pricing problem. These are described in details in \S\ref{pricing}. Other features of the algorithm, such as the branching rules and the generation of cuts, are discussed in \S\ref{branchcut}.

\subsection{Pricing Problem} \label{pricing}
When computing the linear relaxation of the set-partitioning-based formulation from \S\ref{formulations} by column generation, the pricing problem corresponds to finding routes $r\in\mathcal{R}$ (or determining that one does not exist) such that:
\begin{equation} \label{eq:pricing}
c_{r}-\sum_{i=1}^{N}a_{ir}\alpha_{i}-\sum_{k=1}^{C}\sum_{v_{i}\in S_{k}}\sum_{v_{j}\notin S_{k}}b_{r}^{i,j}\beta_{k}<0\,,
\end{equation}
where $\alpha_{i}$, $i\in\{1,\ldots,N\}$, are the dual values associated with the partitioning constraints \eqref{spc}, and $\beta_{k}$, $k\in\{1,\ldots,C\}$, are the dual values associated with the rounded capacity cuts \eqref{rccsp}. If a constraint fixing the number of routes is added to the extended formulation, then its corresponding dual value must also be considered in \eqref{eq:pricing}, and only minor changes in the algorithm are needed. In order to simplify notation we proceed without assuming a fixed number of routes, but we do also consider this case later in the computational experiments.

% HERE
Note that $c_{r}$ is computed by solving a stochastic dynamic program, so it cannot be expressed as a simple function of the individual arcs traversed in $r$. Therefore, we do not attempt to model the pricing problem as a resource-constrained shortest-path problem (RCSPP), as commonly done in other variants of the VRP \citep{irnich2005shortest}. Nevertheless, it is convenient to express the contribution of the dual values to the reduced cost of a route as a function of the individual arcs traversed in the route. To this end, we define $\alpha_{0}=0$, and $\gamma_{ij}=-\alpha_{j}-\sum_{k=1}^{C}\beta_{k}[v_{i}\in S_{k}, v_{j}\notin S_{k}]$. The values $\gamma_{ij}$, $i,j\in\{0,\ldots,N\}$, correspond to the contribution of the dual values to the reduced cost of some route $r$, whenever node $v_{j}$ is visited immediately after $v_{i}$ in $r$.

It is clearly impractical to evaluate $c_{r}$ by solving the full stochastic dynamic program for every route $r$ that could potentially satisfy \eqref{eq:pricing}. Next, we describe an exhaustive search strategy, which partially overcomes the burden of having to calculate (entirely) the dynamic program for every single feasible route.

\subsubsection{Backward Labeling Algorithm.}
The idea of proposing backward labeling for solving the pricing problem is natural, since the dynamic programming recursion is also solved backwardly. We start with $N$ \emph{root} labels, each representing the set of feasible routes \emph{finishing} at some customer. Let $\mathcal{L}$ be a label, and $\mathcal{L}.a$ some attribute $a$ of $\mathcal{L}$. A label possesses various attributes, which are described in Table \ref{tab:attributes}, together with their corresponding initialization values at a root label.

The reduced cost of the route $\{v_{0}\}\oplus\mathcal{L}.pr$, where $\oplus$ denotes sequence concatenation, is denoted by $\phi(\{v_{0}\}\oplus\mathcal{L}.pr)$, and given by:
\begin{multline} \label{eq:redcost}
\phi(\{v_{0}\}\oplus\mathcal{L}.pr)=d_{0\mathcal{L}.n}+\sum_{k=0}^{\infty}\Bigl[(d_{\mathcal{L}.n0}+d_{0\mathcal{L}.n})\Psi(k,Q)\\ +\mathcal{L}.\nu(Q\Psi(k,Q)+Q-k)\Bigr]p_{D_{\mathcal{L}.n}}(k)+\mathcal{L}.\Gamma+\gamma_{0\mathcal{L}.n}\,.
\end{multline}

\begin{table}
\caption{Attributes and initialization values of a label.}
  \centering
  \scalebox{0.8}{
\begin{tabular}{lll}
\hline
Attribute & Description & Initialization$^a$ (root labels) \\
\hline
$n$ & Index of the current customer. & $j$ \\
$rl$ & \makecell[tl]{Remaining load available for \\ further extensions.} & $fQ-E[D_{j}]$ \\
$pr$ & Partial route. & $\{v_{j},v_{0}\}$ \\
$ap$ & A priori cost of $pr$. & $d_{j0}$ \\
$\nu(q)\,,\quad q\in\{0,\ldots,Q\}$ & \makecell[tl]{Expected remaining route cost, \\ considering that customer $v_{n}$ has \\ just been served, and the remaining \\ quantity in the vehicle is $q$.} & $d_{j0}\,,\quad\forall q\in\{0,\ldots,Q\}$ \\
$\Gamma$ & \makecell[tl]{Cumulative contribution of the dual \\ values to the reduced cost of routes \\ created by extending this label.} & $\gamma_{j0}$ \\
\hline
\end{tabular}}
\captionsetup{width=.67\textwidth}
  \caption*{\footnotesize $^a$ Assuming initialization to some customer $v_{j}$.}
  \label{tab:attributes}
\end{table}

For some label $\mathcal{L}$, if $\mathcal{L}.rl\geq E[D_{i}]$ and $v_{i}\notin \mathcal{L}.pr$, then (and only then) label $\mathcal{L}$ can be feasibly extended (backwardly) to $v_{i}$. We denote by $\mathcal{E}$ such extension. The attributes of the new label $\mathcal{E}$ are set as follows:
\begin{align}
\mathcal{E}.n&=i \nonumber \,, \\
\mathcal{E}.rl&=\mathcal{L}.rl-E[D_{i}] \nonumber \,, \\
\mathcal{E}.pr&=\{v_{i}\}\oplus\mathcal{L}.pr \nonumber \,, \\
\mathcal{E}.ap&=\mathcal{L}.ap+d_{i\mathcal{L}.n} \nonumber \,, \\
\mathcal{E}.\Gamma&=\mathcal{L}.\Gamma+\gamma_{i\mathcal{L}.n} \nonumber \,.
\end{align}

We set the values $\mathcal{E}.\nu(q)$, $q\in\{0,\ldots,Q\}$, by performing one step of the stochastic dynamic programming algorithm:
\begin{multline} \label{eq:updatevq}
\mathcal{E}.\nu(q)=\min\Biggl\{d_{i\mathcal{L}.n}+\sum_{k=0}^{\infty}\left[(d_{\mathcal{L}.n0}+d_{0\mathcal{L}.n})\Psi(k,q)+\mathcal{L}.\nu(Q\Psi(k,q)+q-k)\right]p_{D_{\mathcal{L}.n}}(k),\\d_{i0}+d_{0\mathcal{L}.n}+\sum_{k=0}^{\infty}\left[(d_{\mathcal{L}.n0}+d_{0\mathcal{L}.n})\Psi(k,Q)+\mathcal{L}.\nu(Q\Psi(k,Q)+Q-k)\right]p_{D_{\mathcal{L}.n}}(k)\Biggr\} \,.
\end{multline}

Note that, for some $q$, the updating step in \eqref{eq:updatevq} has a time complexity linear on the number of demand values of customer $v_{\mathcal{L}.n}$ that can occur with positive probability. In particular, the recursive step of the dynamic programming is computed in constant time (for a given demand value $k$), since this result is already stored in label $\mathcal{L}$.

The labeling algorithm extends each label (starting with the root labels) to every possible customer. At every new label that is created, the algorithm checks, by evaluating \eqref{eq:redcost}, whether the label leads to a column with a negative reduced cost. We discuss how to control the combinatorial growth of labels in \S\ref{subsecdominance} and \S\ref{subseccbounds}.

% example
By storing in each label the intermediate results of the dynamic program, we avoid solving it entirely for every single route. For example, consider some label $\mathcal{L}$, such that $\mathcal{L}.pr=\{v_{2}, v_{1}, v_{0}\}$. If we want to evaluate the reduced cost of the routes $\{v_{0}, v_{2}, v_{1}, v_{0}\}$ and $\{v_{0}, v_{3}, v_{2}, v_{1}, v_{0}\}$, then we need to evaluate expression \eqref{eq:updatevq} once for every $q\in\{0,\ldots,Q\}$ (backward extension from $v_{2}$ to $v_{3}$), and expression \eqref{eq:redcost} twice (backward extensions from $v_{2}$ and $v_{3}$ to $v_{0}$, completing the routes). If the intermediate results were not stored in the labels, then the computational effort to evaluate the reduced cost of these two routes would be significantly higher, as the full stochastic dynamic program would have to be solved for each of them.

\subsubsection{\emph{ng}-Route Relaxation and Label Dominance.} \label{subsecdominance}
If measures to control the combinatorial growth of labels are not taken, solving the pricing problem becomes a computationally intractable task, even for small VRPSD instances. In order to propose label dominance rules, we relax the route elementarity constraint, i.e., the constraint that enforces a customer to be visited not more than once in a route.

The \emph{ng}-route relaxation \citep{baldacci2011new} is regarded as one of the best route elementarity relaxations for using in branch-and-price algorithms for the VRP \citep{poggi2014new}. It consists in defining for every customer $v_{i}$ a set $\eta_{i}$ containing $v_{i}$ and some other customers (e.g., the 10 nearest customers to $v_{i}$). With \emph{ng}-routes, some label $\mathcal{L}$ can also be (backwardly) extended to a customer $v_{i}\in\mathcal{L}.pr$ (i.e., a customer already visited), but only if some customer $v_{j}$ such that $v_{i}\notin\eta_{j}$ has been visited before the first visit to $v_{i}$ in $\mathcal{L}.pr$ \citep[see][]{baldacci2011new}. An attribute $ng$ is added to the labels, representing the set of customers to which a backward extension is forbidden. The routes that can be generated by these extension rules are simply called \emph{ng}-feasible. Note that every feasible route to the VRPSD is also \emph{ng}-feasible.

When elementarity is relaxed, the coefficients $a_{ir}$ of the partitioning constraints indicate the number of times that customer $v_{i}$ is visited in route $r$, and the coefficients $b_{r}^{i,j}$ of the rounded capacity cuts indicate the number of times that $v_{j}$ is visited immediately after $v_{i}$ in $r$. Note that the partitioning constraints prevent any route visiting the same customer more than once from appearing in a feasible integer solution. We next define label dominance, and then show that dominated labels can be discarded when solving the pricing problem with the labeling algorithm.
\begin{definition} \label{defdominance}
Let $\mathcal{L}^{1}$ and $\mathcal{L}^{2}$ be two labels such that $\mathcal{L}^{1}.n=\mathcal{L}^{2}.n$. We say that $\mathcal{L}^{1}$ \emph{dominates} $\mathcal{L}^{2}$, or alternatively that $\mathcal{L}^{2}$ \emph{is dominated by} $\mathcal{L}^{1}$, if $\mathcal{L}^{1}.rl\geq\mathcal{L}^{2}.rl$, $\mathcal{L}^{1}.ng\subset\mathcal{L}^{2}.ng$, $\mathcal{L}^{1}.\Gamma\leq\mathcal{L}^{2}.\Gamma$, and $\mathcal{L}^{1}.\nu(q)\leq\mathcal{L}^{2}.\nu(q)$, $\forall q\in\{0,\ldots,Q\}$.
\end{definition}

Once two labels $\mathcal{L}^{1}$ and $\mathcal{L}^{2}$ have a dominance relation (i.e., $\mathcal{L}^{1}$ dominates, or is dominated by, $\mathcal{L}^{2}$), further common extensions preserve this relation. We formalize this result with the following proposition:
\begin{proposition} \label{proppreservedom}
Let $\mathcal{L}^{1}$ and $\mathcal{L}^{2}$ be two labels such that $\mathcal{L}^{1}$ dominates $\mathcal{L}^{2}$. Consider an extension $e=\{v_{s_{1}},\ldots,v_{s_{l}}\}$, $l\geq 1$, $s_{i}\in\{1,\ldots,N\}$. Let $\mathcal{E}^{1}$ and $\mathcal{E}^{2}$ be, respectively, the labels corresponding to the extensions of $\mathcal{L}^{1}$ and $\mathcal{L}^{2}$ to $\{v_{s_{1}},\ldots,v_{s_{l}}\}$. Then, $\mathcal{E}^{1}$ dominates $\mathcal{E}^{2}$.
\end{proposition}
\proof
Consider a single-customer extension $e=\{v_{i}\}$. Let $\mathcal{E}^{1}$ and $\mathcal{E}^{2}$ be, respectively, the labels corresponding to the extensions of $\mathcal{L}^{1}$ and $\mathcal{L}^{2}$ to $\{v_{i}\}$. Clearly, $\mathcal{E}^{1}.n=\mathcal{E}^{2}.n=i$, $\mathcal{E}^{1}.rl\geq\mathcal{E}^{2}.rl$, and $\mathcal{E}^{1}.\Gamma\leq\mathcal{E}^{2}.\Gamma$. The values $\mathcal{E}^{1}.\nu(q)$ and $\mathcal{E}^{2}.\nu(q)$, $q\in\{0,\ldots,Q\}$, are defined according to \eqref{eq:updatevq}. When updating $\mathcal{E}^{1}.\nu(q)$ and $\mathcal{E}^{2}.\nu(q)$, the only terms in the minimizer that differ are the expected values:
\begin{align}
\sum_{k=0}^{\infty}&\mathcal{L}.\nu(Q\Psi(k,q)+q-k)p_{D_{\mathcal{L}.n}}(k)\,, \nonumber \\
\intertext{in the first argument to the minimizer, and}
\sum_{k=0}^{\infty}&\mathcal{L}.\nu(Q\Psi(k,Q)+Q-k)p_{D_{\mathcal{L}.n}}(k) \,, \nonumber
\end{align}
in the second argument.

When updating $\mathcal{E}^{1}$ we have $\mathcal{L}=\mathcal{L}^{1}$, and when updating $\mathcal{E}^{2}$ we have $\mathcal{L}=\mathcal{L}^{2}$. Since $\mathcal{L}^{1}.\nu(q)\leq\mathcal{L}^{2}.\nu(q)$, $q\in\{0,\ldots,Q\}$, and also considering nonnegativity of probabilities, it follows that $\mathcal{E}^{1}.\nu(q)\leq\mathcal{E}^{2}.\nu(q)$ must hold, $\forall q\in\{0,\ldots,Q\}$. Therefore, for a single-customer extension $\mathcal{E}^{1}$ dominates $\mathcal{E}^{2}$. Since a multiple-customers extension can be seen as a series of single-customer extensions, the result also holds for that case.
\endproof

The following proposition shows that a dominated label can be safely discarded, without the risk that the labeling algorithm fails to identify a column satisfying \eqref{eq:pricing}, when one exists:
\begin{proposition} \label{propositiondominance}
Let $\mathcal{L}^{1}$ and $\mathcal{L}^{2}$ be two labels such that $\mathcal{L}^{1}$ dominates $\mathcal{L}^{2}$. Let $e=\{v_{s_{1}},\ldots,v_{s_{l}}\}$, $l\geq 0$, $s_{i}\in\{1,\ldots,N\}$, be some (possibly empty) \emph{ng}-feasible extension to $\mathcal{L}^2$. Then, \emph{(i)} route $\{v_{0}\}\oplus e \oplus\mathcal{L}^{1}.pr$ is \emph{ng}-feasible, and \emph{(ii)} its reduced cost is not greater than the reduced cost of route $\{v_{0}\}\oplus e \oplus\mathcal{L}^{2}.pr$.
\end{proposition}
\proof
By the dominance relation we have $\mathcal{L}^{1}.rl\geq\mathcal{L}^{2}.rl$ and $\mathcal{L}^{1}.ng\subset\mathcal{L}^{2}.ng$, which implies that every \emph{ng}-feasible extension to $\mathcal{L}^2$ must also be \emph{ng}-feasible to $\mathcal{L}^1$. Therefore, (i) holds. We consider now (ii). Consider an empty extension $e=\{\}$. By the dominance relation, $\mathcal{L}^{1}.\Gamma\leq\mathcal{L}^{2}.\Gamma$ and $\mathcal{L}^{1}.\nu(q)\leq\mathcal{L}^{2}.\nu(q)$, $\forall q\in\{0,\ldots,Q\}$. Therefore, according to \eqref{eq:redcost}, $\phi(\{v_{0}\}\oplus\mathcal{L}^{1}.pr)\leq\phi(\{v_{0}\}\oplus\mathcal{L}^{2}.pr)$, and (ii) holds for the empty extension. Finally, consider some non-empty extension $e=\{v_{s_{1}},\ldots,v_{s_{l}}\}$, $l\geq 1$. Let $\mathcal{E}^{1}$ and $\mathcal{E}^{2}$ be, respectively, the labels corresponding to the extensions of $\mathcal{L}^{1}$ and $\mathcal{L}^{2}$ to $\{v_{s_{1}},\ldots,v_{s_{l}}\}$. By Proposition \ref{proppreservedom}, $\mathcal{E}^{1}$ dominates $\mathcal{E}^{2}$, so $\mathcal{E}^{1}.\Gamma\leq\mathcal{E}^{2}.\Gamma$ and $\mathcal{E}^{1}.\nu(q)\leq\mathcal{E}^{2}.\nu(q)$, $\forall q\in\{0,\ldots,Q\}$, and thus $\phi(\{v_{0}\}\oplus e\oplus\mathcal{L}^{1}.pr)\leq\phi(\{v_{0}\}\oplus e\oplus\mathcal{L}^{2}.pr)$.
\endproof

\textbf{Heuristic dominance.} The dominance rule from Definition \ref{defdominance} is arguably weak, in the sense that many conditions must hold for dominance to occur. A good heuristic dominance rule is obtained by relaxing the condition $\mathcal{L}^{1}.\nu(q)\leq\mathcal{L}^{2}.\nu(q)$, $\forall q\in\{0,\ldots,Q\}$, and instead requiring only that the expected cost of route $\{v_{0}\}\oplus\mathcal{L}^{1}.pr$ to be less than the expected cost of route $\{v_{0}\}\oplus\mathcal{L}^{2}.pr$. This heuristic rule is useful when solving the pricing problem initially, when there are still many columns satisfying \eqref{eq:pricing}. By applying solely heuristic dominance, it is also possible to transform the whole algorithm into a heuristic. In fact, the results presented in \S\ref{results} show that such heuristic could find the optimal solution in all instances where the optimum is known.

\subsubsection{Reduced cost bounds.} \label{subseccbounds}
Reduced cost bounds, or completion bounds, are lower-bounds computed on the reduced cost of all \emph{ng}-feasible route extensions from a particular label. These bounds are useful for accelerating the solution of the pricing problem, since a label can be safely discarded if it has a nonnegative completion bound. In our algorithm, these bounds play an important role, as in some instances the dominance rules alone are not enough to control the combinatorial growth of labels.

\textbf{Forward RCSP bound.} The expected cost of a \emph{ng}-feasible route is bounded from below by its a priori cost (otherwise, some restocking trip between, say, customers $v_{i}$ and $v_{j}$, would save transportation costs, implying $d_{i0}+d_{0j}<d_{ij}$, which violates the triangle inequality). As a consequence, the expected restocking cost of a route must be nonnegative. The resource-constrained shortest-path (RCSP) bound explores this property, computing a completion bound by underestimating the contribution of the expected route cost to the reduced cost of a column. The following proposition formalizes the RCSP bound:
\begin{proposition}
Let  $\mathcal{G}$ be a graph with node set $\mathcal{N}=\{v_{0},\ldots,v_{N}\}$ and arc set $\mathcal{A}=\{(i,j):v_{i},v_{j}\in\mathcal{N},i\neq j\}$. The cost of each arc $(i,j)\in\mathcal{A}$ is given by $d_{ij}+\gamma_{ij}$. When node $v_{i}$, $i\in\{1,\ldots,N\}$, is visited by some path in $\mathcal{G}$, $E[D_{i}]$ units of the capacity resource are consumed. Let $\psi(i, C)$ be the value of the resource-constrained shortest-path in $\mathcal{G}$ from node $v_{0}$ to node $v_{i}$, assuming a capacity of $C$ when starting the path at $v_{0}$. Then, a lower-bound on the reduced cost of all \emph{ng}-feasible routes that can be generated from some label $\mathcal{L}$ is given by:
\begin{equation} \nonumber
\psi(\mathcal{L}.n,\mathcal{L}.rl+E[D_{\mathcal{L}.n}])+\mathcal{L}.ap+\mathcal{L}.\Gamma\,.
\end{equation}
\end{proposition}
\proof
Let $\mathcal{L}$ be some label, and $e$ the (possibly empty) \emph{ng}-feasible extension to $\mathcal{L}$ that minimizes $\phi(\{v_{0}\}\oplus e\oplus\mathcal{L}.pr)$. Let $\mathcal{E}$ be the label corresponding to the extension of $\mathcal{L}$ to $e$. Assume, for the sake of contradiction, that:
\begin{equation} \label{eq:contrarcspbound}
\phi(r)<\psi(\mathcal{L}.n,\mathcal{L}.rl+E[D_{\mathcal{L}.n}])+\mathcal{L}.ap+\mathcal{L}.\Gamma\,,
\end{equation}
where $r=\{v_{0}\}\oplus\mathcal{E}.pr=\{v_{0}\}\oplus e\oplus\mathcal{L}.pr$. We now decompose $\phi(r)$:
\begin{equation} \nonumber
\phi(r)=\mathcal{Q}+(d_{0\mathcal{E}.n}+\mathcal{E}.ap)+(\gamma_{0\mathcal{E}.n}+\mathcal{E}.\Gamma)\,,
\end{equation}
with $\mathcal{Q}$ representing the expected restocking cost of $r$, the first parenthesis the a priori cost of $r$, and the second parenthesis the total contribution of the dual values to the reduced cost of $r$. Substituting in \eqref{eq:contrarcspbound} and reorganizing the terms, we have:
\begin{equation} \label{eq:reorg}
(d_{0\mathcal{E}.n}+\mathcal{E}.ap-\mathcal{L}.ap)+(\gamma_{0\mathcal{E}.n}+\mathcal{E}.\Gamma-\mathcal{L}.\Gamma)<\psi(\mathcal{L}.n,\mathcal{L}.rl+E[D_{\mathcal{L}.n}])-\mathcal{Q}\,.
\end{equation}

Note that $\{v_{0}\}\oplus e\oplus\{v_{\mathcal{L}.n}\}$ corresponds to a path in $\mathcal{G}$ from $v_{0}$ to $v_{\mathcal{L}.n}$. By \emph{ng}-feasibility of $e$, this path consumes at most $\mathcal{L}.rl+E[D_{\mathcal{L}.n}]$ units of capacity. Moreover, the value of this path is given by the left-hand side of \eqref{eq:reorg}, and (by nonnegativity of $\mathcal{Q}$) it is smaller than $\psi(\mathcal{L}.n,\mathcal{L}.rl+E[D_{\mathcal{L}.n}])$, which is a contradiction. Therefore, assumption \eqref{eq:contrarcspbound} must be wrong.
\endproof

\textbf{Knapsack bound.} This completion bound is computed by solving an unbounded knapsack problem, in which one item is defined for each customer, with a value corresponding to the maximum reduced cost decrease when the associated customer is visited:
\begin{proposition}
Consider an unbounded knapsack instance, where the $N$ items have values $g_{i}=\max\limits_{j\in\{0,\ldots,N\}}(-\gamma_{ji})$, and weights $w_{i}=E[D_{i}]$, $i\in\{1,\ldots,N\}$. Let $\kappa(C)$ be the value of the optimal solution to this knapsack problem, when the capacity of the knapsack is $C$. Then, a lower-bound on the reduced cost of all \emph{ng}-feasible routes that can be generated from some label $\mathcal{L}$ is given by:
\begin{equation} \nonumber
\phi(\{v_{0}\}\oplus\mathcal{L}.pr)-\kappa(\mathcal{L}.rl)\,.
\end{equation}
\end{proposition}
\proof
By the triangle inequality, the expected cost of the route $\{v_{0}\}\oplus\mathcal{L}.pr$ cannot decrease when new customers are added into it. Thus, the value $g_{i}$ is an upper-bound on the reduced cost decrease when adding customer $v_{i}$ to the route $\{v_{0}\}\oplus\mathcal{L}.pr$. Therefore, $\kappa(\mathcal{L}.rl)$ is an upper-bound on the maximum reduced cost decrease that can be achieved when adding customers to route $\{v_{0}\}\oplus\mathcal{L}.pr$.
\endproof

An important characteristic of both the RCSP and the knapsack bounds is that they do not depend on internal information stored in a label, other than its remaining load. This allows them to be pre-processed at the beginning of an iteration of the pricing problem. When the pricing algorithm is called, it uses the updated dual values and pre-computes, for every customer node $\{v_{1},\ldots,v_{N}\}$, and for every possible remaining load $\{0,\ldots,fQ\}$, the values of both completion bounds (i.e., the functions $\psi$ and $\kappa$). This allows access to these bounds in constant time, whenever a new label is created during the labeling algorithm.

\subsection{Branch-and-Bound and Cut Separation} \label{branchcut}
Before every iteration of the pricing problem, we separate possibly violated rounded capacity cuts, and add them to the extended formulation. The labeling algorithm initially enforces only the heuristic dominance from \S\ref{subsecdominance}. Once no more columns can be generated heuristically, we solve the pricing problem exactly by checking exact label dominance. Usually, this step is only required to confirm that no more columns with negative reduced cost exist. In the vast majority of the iterations, the heuristic pricing (i.e., when using heuristic dominance) is able to identify all columns with negative reduced cost.

When the solution to some node of the branch-and-bound tree is fractional and less than the current best upper-bound, we branch on some arc variable of the compact formulation that has a fractional value. We choose the arc variable yielding the highest increase in the lower-bound of the child nodes, when only the columns generated so far are considered. Note that this involves solving up to two (in the worst case) linear programs for every possible candidate arc. These evaluations can be efficiently implemented, by just loading one linear program with all the columns, and changing its cost vector accordingly, i.e., before every solver call, we penalize the columns that do not comply with the branching decisions.

The branch-and-bound tree exploration strategy adopted was that of always choosing the node of the tree with the smallest lower-bound (i.e., best-bound strategy). We have not experimented with other exploration strategies, since we have not observed major inefficiencies when solving the linear relaxation, i.e., there would not be much improvement, in terms of computational time for solving the linear relaxations, if some depth-first strategy \citep[such as the one in][]{desaulniers2014vehicle} were implemented.

Implementation parameters were kept to a minimum: the customer sets used in the \emph{ng}-route relaxation have a fixed size of 12, and contain the own customer plus the 11 nearest customers to it; when solving the pricing problem, we generate a maximum of 20 columns before completing one iteration and retrieving new dual values; finally, a maximum of 4 rounded capacity cuts were added per call to the separation algorithm. These values have been defined after experimenting with different settings in a small number of instances.

\section{Computational Results} \label{results}
In order to assess the effectiveness of the algorithm, we first performed experiments on several literature instances, which had not yet been solved (under optimal restocking) by any existing method. The results are presented in \S\ref{explit}. In \S\ref{louveauxcomp}, we compare results with \cite{louveaux2017}. In \S\ref{dtdcomp}, we compare results with \cite{Gauvinetal2014}, where the detour-to-depot policy was assumed. In \S\ref{loadfactor}, we report the results of the experiments with load factors larger than 1. Finally, in \S\ref{longrun} we describe the results obtained in long-runs of the algorithm on larger instances.

All experiments were performed on a single thread of an Intel\textsuperscript{\textregistered} Xeon\textsuperscript{\textregistered} E5-2670 v2 (2.50GHz) processor, with available memory of 24 gigabytes. Except in the long-runs, the computing time was limited to 5 hours. Linear programs were solved using IBM\textsuperscript{\textregistered} CPLEX\textsuperscript{\textregistered} version 12.6.1. Except for the experiments in \S\ref{louveauxcomp}, where (for comparability) triangular distributions were assumed, the demands of the customers followed Poisson distributions. A precision of $10^{-5}$ was used for the probability computations, i.e., whenever the probability of a certain demand value fell below this threshold, it was truncated to zero.

\subsection{Experiments on literature instances} \label{explit}
We tested the algorithm on instances of the sets $\mathtt{A}$, $\mathtt{E}$, $\mathtt{P}$, and $\mathtt{X}$ of the $\mathtt{CVRPLIB}$ \citep{uchoa2017new}. The demands of the customers were assumed Poisson distributed, with parameters as given by the respective demand values in the instances. Since the computation of expected costs is intrinsically fractional, we did not round the distance between nodes to the nearest integer. We attempted to solve each instance twice: one allowing any number of routes, and the other enforcing the number of routes to be equal to the minimum required, which is equivalent to associating a very high fixed cost for using a vehicle. The results are presented in Table \ref{tablelitAE} (sets $\mathtt{A}$ and $\mathtt{E}$) and Table \ref{tablelitPX} (sets $\mathtt{P}$ and $\mathtt{X}$).

\begin{table}
\caption{Literature instances (sets $\mathtt{A}$ and $\mathtt{E}$), load factor $f=1.00$.}
  \centering
  \scalebox{0.8}{
\begin{tabular}{lrrrrrllrllrl}
\hline
Instance & N$^a$ & K$^b$ & $\frac{N}{K}$ & BP\&C$^c$ & \#$^d$ & gap$^e$ & Time & h-BP\&C$^f$ & Time & VSS$^g$ & K-BP\&C$^h$ & Time$^i$ \\
\hline
A-n32-k5 & 31 & 5 & 6.2 & $^*$856.310 & 5 &  & 0.35 & 856.310 & 0.27 & 4.18\% & = & 0.28 \\
A-n33-k5 & 32 & 5 & 6.4 & $^*$704.801 & 5 &  & 0.01 & 704.801 & 0.01 & 2.43\% & = & 0.01 \\
A-n33-k6 & 32 & 6 & 5.3 & $^*$794.415 & 6 &  & 0.00 & 794.415 & 0.01 & 2.71\% & = & 0.00 \\
A-n34-k5 & 33 & 5 & 6.6 & $^*$828.496 & 6 &  & 0.05 & 828.496 & 0.06 & 1.22\% & $^*$833.777 & 0.02 \\
A-n36-k5 & 35 & 5 & 7.0 & 879.255 & 5 & 2.78\% & 5.00 & 862.309 & 4.18 & 4.91\% & 866.740 & 5.00 \\
A-n37-k5 & 36 & 5 & 7.2 & $^*$710.068 & 5 &  & 2.59 & 710.068 & 0.75 & 0.33\% & = & 2.35 \\
A-n37-k6 & 36 & 6 & 6.0 & $^*$1033.08 & 7 &  & 0.15 & 1033.08 & 0.11 & 3.62\% & $^*$1044.08 & 0.03 \\
A-n38-k5 & 37 & 5 & 7.4 & $^*$779.406 & 6 &  & 0.37 & 779.406 & 0.41 & 6.81\% & $^*$808.874 & 4.21 \\
A-n39-k5 & 38 & 5 & 7.6 & $^*$875.618 & 6 &  & 0.04 & 875.618 & 0.02 & 3.76\% & $^*$887.547 & 1.05 \\
A-n39-k6 & 38 & 6 & 6.3 & $^*$877.922 & 6 &  & 0.11 & 877.922 & 0.21 & 8.56\% & = & 0.13 \\
A-n44-k6 & 43 & 6 & 7.2 & $^*$1026.45 & 7 &  & 4.23 & 1026.45 & 1.89 & 2.15\% & $^*$1029.60 & 0.04 \\
A-n45-k6 & 44 & 6 & 7.3 & $^*$1029.68 & 7 &  & 3.77 & 1029.68 & 1.26 & 5.79\% & 1090.14 & 5.00 \\
A-n45-k7 & 44 & 7 & 6.3 & $^*$1264.94 & 7 &  & 0.77 & 1264.94 & 0.29 & 6.36\% & = & 0.03 \\
A-n46-k7 & 45 & 7 & 6.4 & $^*$1003.23 & 7 &  & 0.27 & 1003.23 & 0.28 & 5.98\% & = & 0.52 \\
A-n48-k7 & 47 & 7 & 6.7 & 1189.10 & 7 & 0.20\% & 5.00 & 1189.10 & 1.25 & 4.81\% & 1198.75 & 5.00 \\
A-n53-k7 & 52 & 7 & 7.4 & $\varnothing$ &  &  & 5.00 & $\varnothing$ & 5.00 &  & $\varnothing$ & 5.00 \\
A-n54-k7 & 53 & 7 & 7.6 & 1292.53 & 8 & 0.12\% & 5.00 & 1292.53 & 5.00 & 3.94\% & 1301.58 & 5.00 \\
A-n55-k9 & 54 & 9 & 6.0 & $^*$1183.93 & 10 &  & 0.29 & 1183.93 & 0.16 & 5.75\% & 1202.70 & 5.00 \\
A-n60-k9 & 59 & 9 & 6.6 & 1536.38 & 10 & 1.58\% & 5.00 & 1536.40 & 5.00 & 4.48\% & 1556.67 & 5.00 \\
A-n61-k9 & 60 & 9 & 6.7 & 1157.86 & 10 & 1.92\% & 5.00 & 1149.00 & 5.00 & 5.68\% & 1199.51 & 5.00 \\
A-n62-k8 & 61 & 8 & 7.6 & $\varnothing$ &  &  & 5.00 & 1442.48 & 5.00 & 6.28\% & $\varnothing$ & 5.00 \\
A-n63-k9 & 62 & 9 & 6.9 & $\varnothing$ &  &  & 5.00 & 1871.12 & 5.00 & 6.32\% & $\varnothing$ & 5.00 \\
A-n63-k10 & 62 & 10 & 6.2 & $\varnothing$ &  &  & 5.00 & 1464.09 & 5.00 & 7.16\% & $\varnothing$ & 5.00 \\
A-n64-k9 & 63 & 9 & 7.0 & $\varnothing$ &  &  & 5.00 & $\varnothing$ & 5.00 &  & $\varnothing$ & 5.00 \\
A-n65-k9 & 64 & 9 & 7.1 & 1318.95 & 10 & 0.36\% & 5.00 & 1318.95 & 5.00 & 7.62\% & 1359.59 & 5.00 \\
A-n69-k9 & 68 & 9 & 7.6 & 1264.42 & 10 & 0.46\% & 5.00 & 1264.42 & 5.00 & 5.83\% & 1282.02 & 5.00 \\
A-n80-k10 & 79 & 10 & 7.9 & $\varnothing$ &  &  & 5.00 & $\varnothing$ & 5.00 &  & $\varnothing$ & 5.00 \\
E-n31-k7 & 30 & 7 & 4.3 & $^*$407.973 & 7 &  & 0.02 & 407.973 & 0.03 & 0.01\% & = & 0.03 \\
E-n51-k5 & 50 & 5 & 10.0 & 557.045 & 5 & 1.23\% & 5.00 & 556.952 & 5.00 & 0.02\% & 556.995 & 5.00 \\
E-n76-k7 & 75 & 7 & 10.7 & 720.872 & 8 & 3.38\% & 5.00 & 703.851 & 5.00 & 3.48\% & 746.534 & 5.00 \\
E-n76-k8 & 75 & 8 & 9.4 & 780.019 & 9 & 1.59\% & 5.00 & 778.509 & 5.00 & 2.31\% & 804.390 & 5.00 \\
E-n76-k10 & 75 & 10 & 7.5 & 894.403 & 11 & 0.56\% & 5.00 & 891.224 & 5.00 & 2.94\% & 918.224 & 5.00 \\
E-n76-k14 & 75 & 14 & 5.4 & $^*$1122.36 & 16 &  & 3.80 & 1122.36 & 3.59 & 5.72\% & 1160.19 & 5.00 \\
E-n101-k8 & 100 & 8 & 12.5 & $\varnothing$ &  &  & 5.00 & $\varnothing$ & 5.00 &  & $\varnothing$ & 5.00 \\
E-n101-k14 & 100 & 14 & 7.1 & 1182.99 & 15 & 0.76\% & 5.00 & 1182.99 & 5.00 & 5.28\% & 1187.79 & 5.00 \\
\hline
\end{tabular}}
\captionsetup{width=.92\textwidth}
  \caption*{\footnotesize $^a$ Number of customers. $^b$ Minimum number of routes for a load factor $f=1.00$. $^c$ Best solution found by the branch-price-and-cut algorithm (number of routes unfixed); $^*$ indicates optimal solutions, and $\varnothing$ indicates the linear bound at the root node could not be computed. $^d$ Number of routes in the best solution of BP\&C. $^e$ Optimality gap of BP\&C. $^f$ Best solution found when only heuristic dominance was enforced (number of routes unfixed). $^g$ Value of the stochastic solution, considering the best solution between BP\&C and h-BP\&C. $^h$ Best solution found by the branch-price-and-cut algorithm (number of routes fixed to $K$); = indicates the solution is the same as in BP\&C. \mbox{$^i$ Running} time in hours.}
  \label{tablelitAE}
\end{table}

\begin{table}
\caption{Literature instances (sets $\mathtt{P}$ and $\mathtt{X}^a$), load factor $f=1.00$.}
  \centering
  \scalebox{0.8}{
\begin{tabular}{lrrrrrllrllrl}
\hline
Instance & N & K & $\frac{N}{K}$ & BP\&C & \# & gap & Time & h-BP\&C & Time & VSS & K-BP\&C & Time \\
\hline
P-n16-k8 & 15 & 8 & 1.9 & $^*$514.640 & 8 &  & 0.00 & 514.640 & 0.00 & 0.00\% & = & 0.00 \\
P-n19-k2 & 18 & 2 & 9.0 & $^*$225.641 & 3 &  & 0.04 & 225.641 & 0.05 & 1.58\% & $^*$229.236 & 0.02 \\
P-n20-k2 & 19 & 2 & 9.5 & $^*$234.342 & 2 &  & 0.22 & 234.342 & 0.17 & 0.01\% & = & 0.04 \\
P-n21-k2 & 20 & 2 & 10.0 & $^*$220.531 & 2 &  & 0.00 & 220.531 & 0.00 & 0.01\% & = & 0.05 \\
P-n22-k2 & 21 & 2 & 10.5 & $^*$232.716 & 2 &  & 0.37 & 232.716 & 0.37 & 0.01\% & = & 0.12 \\
P-n22-k8 & 21 & 8 & 2.6 & $^*$592.339 & 9 &  & 0.00 & 592.339 & 0.00 & 0.02\% & $^*$618.099 & 0.00 \\
P-n23-k8 & 22 & 8 & 2.8 & $^*$621.725 & 9 &  & 0.00 & 621.725 & 0.00 & 5.32\% & $^*$656.184 & 0.00 \\
P-n40-k5 & 39 & 5 & 7.8 & $^*$475.705 & 5 &  & 0.00 & 475.705 & 0.00 & 0.67\% & = & 0.02 \\
P-n45-k5 & 44 & 5 & 8.8 & 537.237 & 5 & 0.32\% & 5.00 & 537.237 & 5.00 & 1.98\% & $^*$537.237 & 4.71 \\
P-n50-k7 & 49 & 7 & 7.0 & $^*$587.314 & 7 &  & 0.11 & 587.314 & 0.47 & 3.95\% & = & 0.09 \\
P-n50-k8 & 49 & 8 & 6.1 & $^*$673.430 & 9 &  & 1.49 & 673.430 & 0.68 & 5.03\% & 713.671 & 5.00 \\
P-n50-k10 & 49 & 10 & 4.9 & $^*$763.212 & 11 &  & 3.48 & 763.212 & 3.11 & 3.87\% & $^*$775.579 & 0.08 \\
P-n51-k10 & 50 & 10 & 5.0 & $^*$814.687 & 11 &  & 0.02 & 814.687 & 0.01 & 5.02\% & 850.742 & 5.00 \\
P-n55-k7 & 54 & 7 & 7.7 & $^*$591.846 & 7 &  & 1.64 & 591.846 & 0.59 & 4.34\% & = & 0.70 \\
P-n55-k10 & 54 & 10 & 5.4 & $^*$746.080 & 10 &  & 0.84 & 746.080 & 0.61 & 6.69\% & = & 0.77 \\
P-n55-k15 & 54 & 15 & 3.6 & $^*$1073.11 & 18 &  & 0.01 & 1073.11 & 0.01 & 4.27\% & 1176.60 & 5.00 \\
P-n60-k10 & 59 & 10 & 5.9 & 805.176 & 11 & 0.06\% & 5.00 & 805.176 & 3.96 & 3.70\% & $^*$814.436 & 0.28 \\
P-n60-k15 & 59 & 15 & 3.9 & $^*$1089.23 & 16 &  & 0.02 & 1089.23 & 0.01 & 4.00\% & $^*$1109.88 & 0.18 \\
P-n65-k10 & 64 & 10 & 6.4 & 859.843 & 10 & 0.69\% & 5.00 & 858.768 & 5.00 & 0.87\% & $^*$857.481 & 1.19 \\
P-n70-k10 & 69 & 10 & 6.9 & 888.060 & 11 & 0.22\% & 5.00 & 888.352 & 5.00 & 4.82\% & 926.702 & 5.00 \\
P-n76-k4 & 75 & 4 & 18.8 & $\varnothing$ &  &  & 5.00 & $\varnothing$ & 5.00 &  & $\varnothing$ & 5.00 \\
P-n76-k5 & 75 & 5 & 15.0 & $\varnothing$ &  &  & 5.00 & $\varnothing$ & 5.00 &  & $\varnothing$ & 5.00 \\
P-n101-k4 & 100 & 4 & 25.0 & $\varnothing$ &  &  & 5.00 & $\varnothing$ & 5.00 &  & $\varnothing$ & 5.00 \\
X-n101-k25 & 100 & 25 & 4.0 & 30685.8 & 29 & 0.79\% & 5.00 & 30738.2 & 5.00 & 9.01\% & 37452.4 & 5.00 \\
X-n110-k13 & 109 & 13 & 8.4 & 16760.8 & 14 & 0.57\% & 5.00 & 16734.1 & 5.00 & 4.16\% & 17201.1 & 5.00 \\
X-n148-k46 & 147 & 46 & 3.2 & 56584.4 & 50 & 0.03\% & 5.00 & 56584.4 & 2.14 & 1.69\% & 57223.5 & 5.00 \\
X-n195-k51 & 194 & 51 & 3.8 & 50067.1 & 58 & 0.66\% & 5.00 & 50025.2 & 5.00 & 10.48\% & $\varnothing$ & 5.00 \\
\hline
\end{tabular}}
\captionsetup{width=.92\textwidth}
  \caption*{\footnotesize $^a$ All instances with up to 200 nodes were tested. Results are only reported for the instances where at least one feasible solution could be found.}
  \label{tablelitPX}
\end{table}

Overall, the algorithm was effective in most instances with up to 60 nodes. The largest instance solved to optimality was $\mathtt{E\text{-}n76\text{-}k14}$, with 76 nodes. Among the unsolved instances, the optimality gap was less than 1\% in 14 of them, including 4 instances of the set $\mathtt{X}$ with 101 to 195 nodes. In some instances (for example, $\mathtt{A\text{-}n53\text{-}k7}$), the algorithm failed to compute the linear bound at the root node of the branch-and-bound tree. In these few cases, the dominance rules and bounding procedures were ineffective to control the combinatorial growth of labels.

The value of the stochastic solution (VSS) was computed by choosing the best solution between BP\&C and h-BP\&C. The deterministic equivalent solution corresponds to the optimal capacitated VRP solution, when the number of routes is not fixed. Each route of the deterministic equivalent solution was evaluated by the optimal restocking algorithm (stochastic dynamic program) in its both possible orientations, and the best one was considered. Note that the two highest VSS values (in excess of 9\%) were obtained in instances with a customers per route ratio below or equal to 4.0. The average VSS over all instances where a feasible solution could be found was 4.0\%. Therefore, taking stochasticity into account when solving these literature instances may lead to significantly better solutions.

Apparently, opening new routes seems to be a good strategy in many cases: from the 32 instances solved to optimality, 16 have optimal solutions with more routes than the minimum required. Opening new routes reduces the average load over all the routes, and decreases the chances of replenishments. As a result, the uncertainty in the solution cost decreases, since the deterministic a priori cost becomes much more relevant than the stochastic restocking cost.

Among the instances solved to optimality in both cases of fixed and unfixed number of routes, 10 of them have different optimal solutions, i.e., the number of routes in the optimal solutions differ. We further investigate these instances in Table \ref{tableunfixedvsfixed}. The highest percentage variation in the solution value occurred in instance $\mathtt{P\text{-}n23\text{-}k8}$, where the solution with the minimum number of routes is 5.5\% worse than the solution when the number of routes is not fixed. We confirm that opening new routes, in fact, reduces the uncertainty in the solutions: the expected restocking cost almost halves when one extra route is allowed. Therefore, opening new routes is an effective insurance against stochastic restocking costs, and in many cases allows superior solutions to the problem.

\begin{table}
\caption{Comparison of optimal solutions for fixed and unfixed number of routes.}
  \centering
  \scalebox{0.8}{
\begin{tabular}{llllllllll}
\hline
Instance & N & K & $\frac{N}{K}$ & BC\&P & $E[\mathcal{Q}]$$^a$ & \# & K-BC\&P & $E[\mathcal{Q}]$$^b$ & \% Diff$^c$ \\
\hline
A-n34-k5 & 33 & 5 & 6.6 & 828.496 & 4.5\% & 6 & 833.777 & 5.0\% & 0.6\% \\
A-n37-k6 & 36 & 6 & 6.0 & 1033.08 & 2.7\% & 7 & 1044.08 & 5.7\% & 1.1\% \\
A-n38-k5 & 37 & 5 & 7.4 & 779.406 & 1.4\% & 6 & 808.874 & 8.0\% & 3.8\% \\
A-n39-k5 & 38 & 5 & 7.6 & 875.618 & 4.2\% & 6 & 887.547 & 5.7\% & 1.4\% \\
A-n44-k6 & 43 & 6 & 7.2 & 1026.45 & 5.2\% & 7 & 1029.60 & 8.8\% & 0.3\% \\
P-n19-k2 & 18 & 2 & 9.0 & 225.641 & 2.2\% & 3 & 229.236 & 7.2\% & 1.6\% \\
P-n22-k8 & 21 & 8 & 2.6 & 592.339 & 0.5\% & 9 & 618.099 & 2.4\% & 4.3\% \\
P-n23-k8 & 22 & 8 & 2.8 & 621.725 & 10.8\% & 9 & 656.184 & 19.0\% & 5.5\% \\
P-n50-k10 & 49 & 10 & 4.9 & 763.212 & 4.0\% & 11 & 775.579 & 8.4\% & 1.6\% \\
P-n60-k15 & 59 & 15 & 3.9 & 1089.23 & 6.4\% & 16 & 1109.88 & 10.3\% & 1.9\% \\
\textbf{Average} &  &  &  &  & \textbf{4.2\%} &  &  & \textbf{8.1\%} & \textbf{2.2\%} \\
\hline
\end{tabular}}
\captionsetup{width=.66\textwidth}
  \caption*{\footnotesize $^a$ Expected restocking cost as a percentage of total solution cost (number of routes not fixed). $^b$ Expected restocking cost as a percentage of total solution cost (number of routes fixed to $K$). \mbox{$^c$ Percentage} increase in the solution value, when the number of routes is fixed.}
  \label{tableunfixedvsfixed}
\end{table}

Finally, when branch-price-and-cut operated in heuristic mode (by enforcing only heuristic dominance), the optimal solutions of all 32 instances solved to optimality could be found. In many cases, the computing time was significantly lower. In other few cases, however, the heuristic took additional time when compared to the exact algorithm (for example, instance $\mathtt{P\text{-}n50\text{-}k7}$). This is explained by different choices of the heuristic and exact algorithms regarding the branching decisions.

\subsection{Comparison with \cite{louveaux2017}} \label{louveauxcomp}
The instances proposed by \cite{louveaux2017} exhibit a high ratio of average number of customers per route, a characteristic that favors pure branch-and-cut solution methods \citep{poggi2014new}. We tested the branch-price-and-cut algorithm on the $32$ instances with up to $51$ nodes. Our algorithm is not effective for the remaining instances, since the ratio of customers per route is just too high. Differently from \cite{louveaux2017}, we do not use bounding procedures tailored for the case of identically distributed demands. Still, the branch-cut-and-price algorithm was able to close three instances that could not be solved by the integer L-shaped algorithm, improving the best known solution in two of them. We also computed the deterministic equivalent solution and measured the (previously unknown) VSS of all 32 instances tested. The VSS has been computed by considering the best known solution for each instance. The results are presented in Table \ref{tablecomplouveaux}.

\begin{table}
\caption{Comparison with \cite{louveaux2017}.}
  \centering
  \scalebox{0.8}{
\begin{tabular}{lrrrrcrrllrl}
\hline
Instance & N & K & N/K & Q & Dist$^a$ & LS-17$^b$ & BP\&C & Time & VSS$^c$ & h-BP\&C & Time \\
\hline
E031-09h & 30 & 2 & 15.0 & 84 & 3 & $^*$332.753 & $^*$332.753 & 0.05 & 0.02\% & 332.753 & 0.09 \\
E031-09h & 30 & 2 & 15.0 & 79 & 3 & $^*$335.296 & $^*$335.296 & 0.07 & 0.00\% & 335.296 & 0.10 \\
E031-09h & 30 & 2 & 15.0 & 84 & 9 & $^*$337.674 & $^*$337.674 & 0.56 & 0.16\% & 337.674 & 0.17 \\
E031-09h & 30 & 2 & 15.0 & 79 & 9 & $^*$344.525 & $^*$344.525 & 0.74 & 0.00\% & 344.525 & 0.52 \\
E031-09h & 30 & 3 & 10.0 & 59 & 3 & $^*$358.947 & $^*$358.947 & 0.20 & 0.00\% & 358.947 & 0.13 \\
E031-09h & 30 & 3 & 10.0 & 56 & 3 & $^*$364.066 & $^*$364.066 & 1.44 & 0.75\% & 364.066 & 1.45 \\
E031-09h & 30 & 3 & 10.0 & 59 & 9 & $^*$367.155 & $^*$367.155 & 0.35 & 0.16\% & 367.155 & 0.16 \\
E031-09h & 30 & 3 & 10.0 & 56 & 9 & 373.129 & $^*$372.784 & 0.91 & 0.60\% & 372.784 & 0.19 \\
\hline
E051-05e & 50 & 2 & 25.0 & 139 & 3 & $^*$441.000 & $\varnothing$ & 5.00 & 0.09\% & $\varnothing$ & 5.00 \\
E051-05e & 50 & 2 & 25.0 & 132 & 3 & $^*$441.311 & $\varnothing$ & 5.00 & 0.00\% & 456.336 & 5.00 \\
E051-05e & 50 & 2 & 25.0 & 139 & 9 & $^*$443.006 & $\varnothing$ & 5.00 & 0.35\% & $\varnothing$ & 5.00 \\
E051-05e & 50 & 2 & 25.0 & 132 & 9 & 448.083 & $\varnothing$ & 5.00 & 0.00\% & 469.676 & 5.00 \\
E051-05e & 50 & 3 & 16.7 & 99 & 3 & $^*$459.000 & 483.009 & 5.00 & 0.24\% & 459.000 & 5.00 \\
E051-05e & 50 & 3 & 16.7 & 93 & 3 & $^*$459.049 & 459.049 & 5.00 & 0.11\% & 459.049 & 2.83 \\
E051-05e & 50 & 3 & 16.7 & 99 & 9 & $^*$460.550 & 482.768 & 5.00 & 1.25\% & 460.550 & 2.58 \\
E051-05e & 50 & 3 & 16.7 & 93 & 9 & $^*$465.629 & 487.260 & 5.00 & 0.25\% & 465.629 & 2.01 \\
\hline
A034-02f & 33 & 2 & 16.5 & 92 & 3 & $^*$1404.64 & 1404.64 & 5.00 & 0.81\% & 1404.64 & 1.80 \\
A034-02f & 33 & 2 & 16.5 & 87 & 3 & $^*$1418.72 & $^*$1418.72 & 3.23 & 0.90\% & 1418.72 & 0.76 \\
A034-02f & 33 & 2 & 16.5 & 92 & 9 & $^*$1436.02 & 1466.40 & 5.00 & 1.24\% & 1436.02 & 1.15 \\
A034-02f & 33 & 2 & 16.5 & 87 & 9 & 1484.13 & 1546.21 & 5.00 & 0.49\% & 1484.13 & 5.00 \\
A034-02f & 33 & 3 & 11.0 & 65 & 3 & $^*$1557.84 & $^*$1557.84 & 0.14 & 4.37\% & 1557.84 & 0.03 \\
A034-02f & 33 & 3 & 11.0 & 62 & 3 & $^*$1595.84 & $^*$1595.84 & 0.10 & 1.41\% & 1595.84 & 0.04 \\
A034-02f & 33 & 3 & 11.0 & 65 & 9 & 1613.76 & $^*$1613.76 & 0.30 & 4.86\% & 1613.76 & 0.10 \\
A034-02f & 33 & 3 & 11.0 & 62 & 9 & 1679.79 & $^*$1648.70 & 0.10 & 4.07\% & 1648.70 & 0.07 \\
\hline
A048-03f & 47 & 2 & 23.5 & 131 & 3 & $^*$1812.02 & $\varnothing$ & 5.00 & 0.00\% & $\varnothing$ & 5.00 \\
A048-03f & 47 & 2 & 23.5 & 124 & 3 & $^*$1818.16 & $\varnothing$ & 5.00 & 0.48\% & $\varnothing$ & 5.00 \\
A048-03f & 47 & 2 & 23.5 & 131 & 9 & $^*$1827.96 & $\varnothing$ & 5.00 & 0.41\% & $\varnothing$ & 5.00 \\
A048-03f & 47 & 2 & 23.5 & 124 & 9 & $^*$1864.33 & $\varnothing$ & 5.00 & 0.75\% & $\varnothing$ & 5.00 \\
A048-03f & 47 & 3 & 15.7 & 93 & 3 & $^*$1953.15 & $\varnothing$ & 5.00 & 0.45\% & 1953.15 & 5.00 \\
A048-03f & 47 & 3 & 15.7 & 88 & 3 & $^*$1960.17 & 2045.11 & 5.00 & 0.22\% & 1960.17 & 1.54 \\
A048-03f & 47 & 3 & 15.7 & 93 & 9 & $^*$1968.48 & $\varnothing$ & 5.00 & 2.72\% & 1968.48 & 4.33 \\
A048-03f & 47 & 3 & 15.7 & 88 & 9 & 1989.60 & 2082.20 & 5.00 & 0.66\% & 1989.60 & 1.74 \\
\hline
\end{tabular}}
\captionsetup{width=.82\textwidth}
  \caption*{\footnotesize $^a$ Number of states of the triangular distribution assumed for the customer demands (all customers i.i.d.). $^b$ Best solution found by the integer L-shaped algorithm from \cite{louveaux2017}. $^c$ Value of the stochastic solution, considering the best solution between LS-17 and BP\&C.}
  \label{tablecomplouveaux}
\end{table}

The highest improvements over the deterministic equivalent solutions were found on the instances with an average of 11 customers per route. Opportunely, on instances with up to that average ratio the branch-price-and-cut algorithm performed better (all the eight instances solved) than the integer L-shaped method (five instances solved). Note that in \cite{louveaux2017} the maximum allowed computational time was also 5 hours, but the processor used in their experiments is approximately 48\% faster (when only accounting for clock speed) than the one used by us.

On the other hand, the 24 remaining instances have a VSS significantly low: maximum of 2.72\% (obtained on instance $\mathtt{A048\text{-}03f}$), and average of only 0.48\% (for comparison, the average VSS on the 8 instances with a ratio up to 11 was 2.03\%). There are two possible explanations for this result: first, when there are many customers per route, the optimal policy has ample opportunity to choose a preventive replenishment point with a low marginal restocking cost (i.e., between consecutive customers $v_{i}$ and $v_{j}$ such that $d_{i0}+d_{0j}-d_{ij}$ is relatively small). Second, when there are many customers in a route, the a priori cost simply outweighs the restocking cost, making the optimization of the deterministic cost a good strategy. In these instances, the integer L-shaped method outperformed the branch-price-and-cut algorithm. However, at least under the demand distributions considered, it is clear that finding the optimal solutions in these instances offers a very small payoff.

The only other approach in the literature that considers multiple vehicles and optimal restocking is the one from \cite{salavati2017exact}. In this paper, the same instances from \cite{louveaux2017} are solved, with slightly worse results in terms of computing time and number of instances solved. In addition, randomly generated non-literature instances are also solved. In these instances, the customer demands follow discrete distributions with up to 5 possible values, and the number of customers per route is also relatively high (between 10 and 30). The results were reported in an aggregated form over 880 runs of the algorithm. The individual parameters and solution of each generated instance are not available, so comparison of results (in a per-instance basis) was not possible.

\subsection{Comparison with the detour-to-depot policy} \label{dtdcomp}
In \cite{Gauvinetal2014}, several literature instances were solved assuming the detour-to-depot policy, maximum load factor of 1, unfixed number of routes, and Poisson distributed demands. By comparing the optimal solution values, we were able to quantify the benefit of adopting the optimal restocking policy. The results are presented in Table \ref{tablecompgauvin}. Note that the distances in \cite{Gauvinetal2014} were rounded to the nearest integer. We repeated the experiments with this new setting, and for this reason the optimal values differ slightly from the ones in \S\ref{explit}.

\begin{table}
\caption{Comparison with \cite{Gauvinetal2014}.}
  \centering
  \scalebox{0.8}{
\begin{tabular}{lrrrrrrlrcc}
\hline
Instance & N & K & Q & $\sum\mu_{i}$$^a$ & GDG-14$^b$ & BP\&C & \% Diff$^c$ & \#$^d$ & $\bar{l}(K)$$^e$ & $\bar{l}(\#)$$^f$ \\
\hline
A-n32-k5 & 31 & 5 & 100 & 410 & 853.6 & 853.5 & 0.01\% & 5 & 0.82 & 0.82 \\
A-n33-k5 & 32 & 5 & 100 & 446 & 704.2 & 704.1 & 0.02\% & 5 & 0.89 & 0.89 \\
A-n33-k6 & 32 & 6 & 100 & 541 & 793.9 & 793.5 & 0.05\% & 6 & 0.90 & 0.90 \\
A-n34-k5 & 33 & 5 & 100 & 460 & 826.9 & 826.8 & 0.01\% & \textbf{6} & \textbf{0.92} & \textbf{0.77} \\
A-n37-k5 & 36 & 5 & 100 & 407 & 708.3 & 707.2 & 0.16\% & 5 & 0.81 & 0.81 \\
A-n37-k6 & 36 & 6 & 100 & 570 & 1030.7 & 1030.5 & 0.02\% & \textbf{7} & \textbf{0.95} & \textbf{0.81} \\
A-n38-k5 & 37 & 5 & 100 & 481 & 775.1 & 775.1 & 0.00\% & \textbf{6} & \textbf{0.96} & \textbf{0.80} \\
A-n39-k5 & 38 & 5 & 100 & 475 & 869.2 & 868.8 & 0.05\% & \textbf{6} & \textbf{0.95} & \textbf{0.79} \\
A-n39-k6 & 38 & 6 & 100 & 526 & 876.6 & 876.4 & 0.03\% & 6 & 0.88 & 0.88 \\
A-n44-k6 & 43 & 6 & 100 & 570 & 1025.5 & 1024.3 & 0.12\% & \textbf{7} & \textbf{0.95} & \textbf{0.81} \\
A-n45-k6 & 44 & 6 & 100 & 593 & 1026.7 & 1026.5 & 0.02\% & \textbf{7} & \textbf{0.99} & \textbf{0.85} \\
A-n45-k7 & 44 & 7 & 100 & 634 & 1264.8 & 1264.1 & 0.05\% & 7 & 0.91 & 0.91 \\
A-n46-k7 & 45 & 7 & 100 & 603 & 1002.2 & 1001.9 & 0.03\% & 7 & 0.86 & 0.86 \\
A-n55-k9 & 54 & 9 & 100 & 839 & 1179.1 & 1179.0 & 0.01\% & \textbf{10} & \textbf{0.93} & \textbf{0.84} \\
P-n16-k8 & 15 & 8 & 35 & 246 & 512.8 & 512.8 & 0.01\% & 8 & 0.88 & 0.88 \\
P-n19-k2 & 18 & 2 & 160 & 310 & 224.1 & 224.0 & 0.03\% & \textbf{3} & \textbf{0.97} & \textbf{0.65} \\
P-n20-k2 & 19 & 2 & 160 & 310 & 233.1 & 233.0 & 0.05\% & 2 & 0.97 & 0.97 \\
P-n21-k2 & 20 & 2 & 160 & 298 & 219.0 & 218.8 & 0.08\% & 2 & 0.93 & 0.93 \\
P-n22-k2 & 21 & 2 & 160 & 308 & 231.3 & 230.8 & 0.23\% & 2 & 0.96 & 0.96 \\
P-n23-k8 & 22 & 8 & 40 & 313 & 619.5 & 619.5 & 0.00\% & \textbf{9} & \textbf{0.98} & \textbf{0.87} \\
P-n40-k5 & 39 & 5 & 140 & 618 & 472.5 & 472.3 & 0.05\% & 5 & 0.88 & 0.88 \\
P-n50-k7 & 49 & 7 & 150 & 951 & 582.4 & 582.2 & 0.03\% & 7 & 0.91 & 0.91 \\
P-n50-k8 & 49 & 8 & 120 & 951 & 669.2 & 669.1 & 0.02\% & \textbf{9} & \textbf{0.99} & \textbf{0.88} \\
P-n50-k10 & 49 & 10 & 100 & 951 & 758.8 & 758.6 & 0.02\% & \textbf{11} & \textbf{0.95} & \textbf{0.86} \\
P-n51-k10 & 50 & 10 & 80 & 777 & 809.7 & 809.3 & 0.04\% & \textbf{11} & \textbf{0.97} & \textbf{0.88} \\
P-n55-k7 & 54 & 7 & 170 & 1042 & 588.6 & 588.5 & 0.02\% & 7 & 0.88 & 0.88 \\
P-n55-k10 & 54 & 10 & 115 & 1042 & 742.4 & 742.0 & 0.05\% & 10 & 0.91 & 0.91 \\
P-n55-k15 & 54 & 15 & 70 & 1042 & 1068.1 & 1067.8 & 0.02\% & \textbf{18} & \textbf{0.99} & \textbf{0.83} \\
P-n60-k15 & 59 & 15 & 80 & 1134 & 1085.5 & 1084.9 & 0.05\% & \textbf{16} & \textbf{0.95} & \textbf{0.89} \\
\hline
\end{tabular}}
\captionsetup{width=.74\textwidth}
  \caption*{\footnotesize $^a$ Sum of the expected demand over all customers ($\mu_{i}=E[D_{i}]$). $^b$ Optimal solution from \cite{Gauvinetal2014} (number of routes unfixed). $^c$ Percentage improvement in the solution value when adopting the optimal policy. $^d$ Number of routes in the optimal solution (GDG-14 and BP\&C always have the same number of routes). \textbf{Bold} entries indicate the value is different from $K$. $^e$ Average route relative load (route load divided by $Q$) if the minimum number of routes were enforced. $^f$ Average route relative load in the optimal solution.}
  \label{tablecompgauvin}
\end{table}

The results indicate that (under the above-mentioned assumptions) the benefit of adopting the optimal policy in place of the detour-to-depot policy is close to none. Note that except for three instances of the set $\mathtt{P}$, whenever the average route relative load (route load divided by $Q$) when using the minimum number of routes would amount to 0.92 or more, it became cost-effective to open additional routes, bringing the average relative load to lower levels. This is valid for both policies, since the number of routes in the optimal solutions of all 29 instances is the same. This takes us to a situation identified in \cite{florio2017} for the single-vehicle case: on routes with a small load relative to the vehicle capacity, the probability of having to perform one replenishment trip is relatively small. Therefore, the deterministic a priori cost supersedes the restocking cost, making the choice of which policy to use almost irrelevant. In other words, the advantages of optimal restocking never materialize when the number of routes is not fixed.

We do not expect these results to sustain when the number of routes is fixed to the minimum. In these cases, the route loads may approach $Q$, or even exceed $Q$, if the load factor parameter allows. In these scenarios, the restocking cost becomes significant, and so also the decision of which policy to apply. Unfortunately, there are no results in the literature to compare both policies for this case.

\subsection{Load factors larger than 1} \label{loadfactor}
Relaxations of the load constraint, even if slight, have the effect of making the problem more difficult. This is mainly caused by two side effects: first, a (potential) increase in the customers per route ratio; second, loss of quality of the forward RCSP bound, since the restocking cost becomes more relevant, and the bound given by the a priori cost becomes weaker. Therefore, in these tests we enforced only heuristic dominance, and experimented only with the instances that could be solved relatively fast for a load factor of 1. For each instance, we experimented with 4 additional load factors ranging from $1.05$ to $1.50$. The results are presented in Table \ref{tableloads}.

\begin{table}
\caption{Load factors from $f=1.05$ up to $f=1.50$.}
  \centering
  \scalebox{0.8}{
\begin{tabular}{lrrrrrrrrrr}
\hline
Instance & BP\&C & \# & f=1.05 & \# & f=1.20 & \# & f=1.35 & \# & f=1.50 & \# \\
\hline
A-n32-k5 & 856.310 & 5 &  &  &  &  &  &  & 831.076 & 4 \\
A-n33-k5 & 704.801 & 5 \\
A-n33-k6 & 794.415 & 6 \\
A-n34-k5 & 828.496 & 6 & 823.241 & 5 &  &  &  &  & 816.479 & 4 \\
A-n37-k6 & 1033.08 & 7 &  &  &  &  & 1033.04 & 6 & 1030.52 & 6 \\
A-n38-k5 & 779.406 & 6 \\
A-n39-k5 & 875.618 & 6 &  &  & 874.863 & 5 \\
A-n39-k6 & 877.922 & 6 \\
A-n45-k7 & 1264.94 & 7 \\
A-n46-k7 & 1003.23 & 7 \\
A-n55-k9 & 1183.93 & 10 \\
E-n31-k7 & 407.973 & 7 \\
P-n16-k8 & 514.640 & 8 &  &  &  &  & 513.160 & 8 & 503.906 & 6 \\
P-n19-k2 & 225.641 & 3 & 213.845 & 2 \\
P-n20-k2 & 234.342 & 2 & 227.738 & 2 \\
P-n21-k2 & 220.531 & 2 \\
P-n22-k2 & 232.716 & 2 \\
P-n22-k8 & 592.339 & 9 &  &  &  &  & 592.317 & 8 & 591.793 & 8 \\
P-n23-k8 & 621.725 & 9 & 620.634 & 9 &  &  & 620.426 & 9 & 613.591 & 7 \\
P-n40-k5 & 475.705 & 5 \\
P-n50-k7 & 587.314 & 7 \\
P-n51-k10 & 814.687 & 11 &  &  &  &  &  &  & 814.569 & 10 \\
P-n55-k15 & 1073.11 & 18 &  &  &  &  & 1069.63 & 17 & 1066.90 & 15 \\
P-n60-k15 & 1089.23 & 16 &  &  &  &  &  &  & 1089.07 & 15 \\
\hline
\end{tabular}}
\captionsetup{width=.72\textwidth}
  \caption*{\footnotesize \textbf{Note:} Only solutions with a strict improvement over previous ones (smaller load factors or BP\&C) are reported.}
  \label{tableloads}
\end{table}

Small increases in the load factor (relative to the base value of 1) may lead to strictly better solutions, with lower total expected cost, fewer routes, or both. By allowing a load factor of $1.50$, the optimal solutions of 12 out of the 24 instances experimented could be improved. In practical scenarios where the maximum route duration is constrained (e.g., due to working hours), it is unlikely that the load factor can be increased arbitrarily. In these situations, it might be more reasonable to enforce only constraints on the maximum expected route cost (or duration), and drop capacity-based constraints altogether.

\begin{figure}
\centering
  \includegraphics[trim={0 200 0 20},clip,width=450pt]{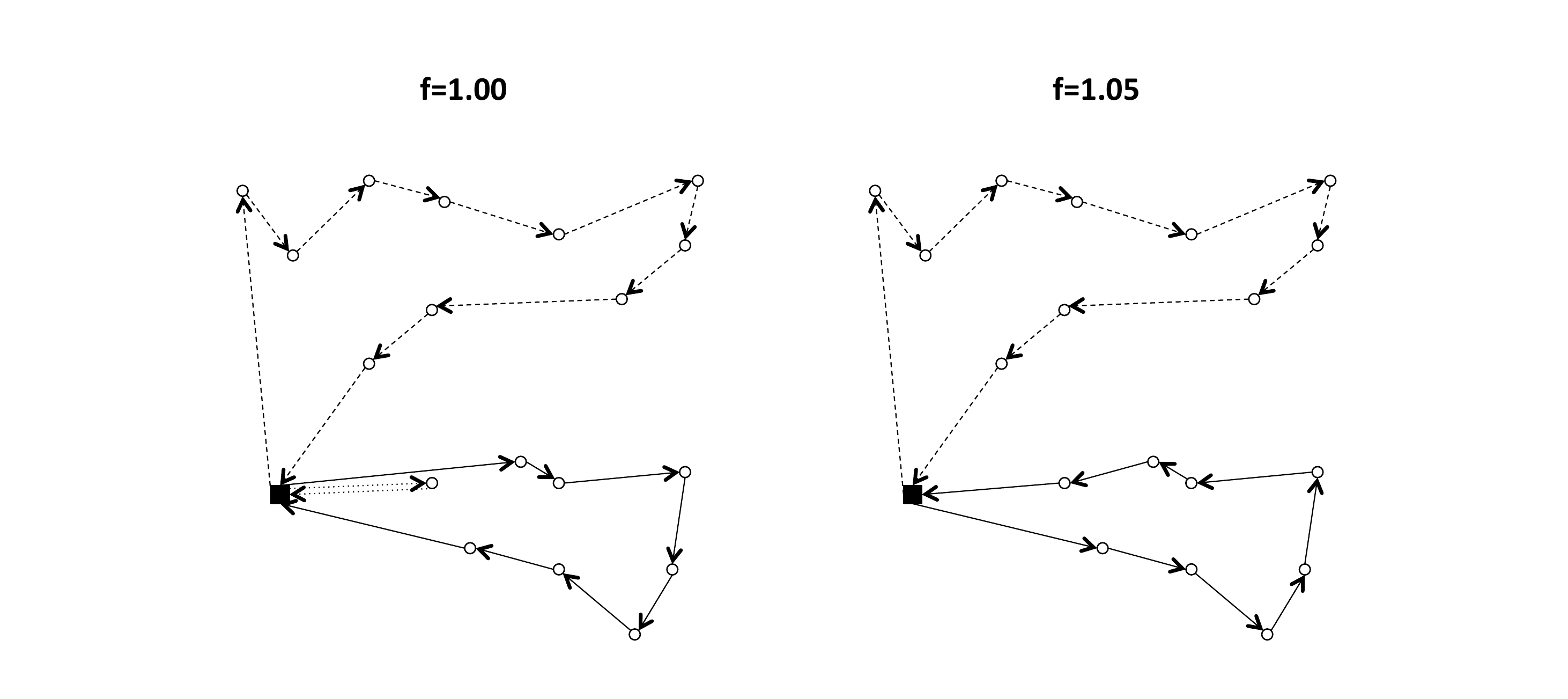}
  \caption{Optimal solutions of instance $\mathtt{P\text{-}n19\text{-}k2}$ when $f\text{=}1.00$ and when $f\text{=}1.05$.}
  \label{figpn19}
\end{figure}

Figure \ref{figpn19} illustrates the effect of allowing a slightly larger load factor in instance $\mathtt{P\text{-}n19\text{-}k2}$. The solution when $f\text{=}1.05$ (confirmed to be optimal) has a value 5.2\% lower and requires one less vehicle than the optimal solution when $f\text{=}1.00$. When $f\text{=}1.05$, the expected cost of the routes are 121.182 (dashed arrows) and 92.663 (solid arrows). Therefore, the solution when $f\text{=}1.05$ is superior regarding total expected cost and number of vehicles, and is not inferior regarding maximum expected route cost (which, in both cases, is 121.182).

\subsection{Long-run experiments} \label{longrun}
In the long-run experiments, we selected the instances where the optimality gap (considering unfixed number of routes) in the first run was below 1\%. We ran again the branch-price-and-cut algorithm on these instances, but this time allowing a maximum running time of 100 hours. The results are presented in Table \ref{tablelongrun}.

\begin{table}
\caption{Long-run results.}
  \centering
  \scalebox{0.8}{
\begin{tabular}{lllllllllllll}
\hline
Instance & N & K & N/K & BP\&C & \# & gap & Time & VSS$^a$ & K-BP\&C & \# & gap & Time \\
\hline
A-n48-k7 & 47 & 7 & 6.7 & 1189.10$^*$ & 7 &  & 12.83 & 4.81\% & 1189.10$^*$ & 7 &  & 19.16 \\
A-n54-k7 & 53 & 7 & 7.6 & 1292.53$^*$ & 8 &  & 11.23 & 3.94\% & 1301.58$^*$ & 7 &  & 6.84 \\
A-n65-k9 & 64 & 9 & 7.1 & 1318.95$^*$ & 10 &  & 64.03 & 7.62\% & 1356.85$^*$ & 9 &  & 42.62 \\
A-n69-k9 & 68 & 9 & 7.6 & 1264.42$^*$ & 10 &  & 80.49 & 5.83\% & 1280.59 & 9 & 0.89\% & 100.00 \\
E-n76-k10 & 75 & 10 & 7.5 & 891.224$^*$ & 11 &  & 70.66 & 2.94\% & 912.839 & 10 & 0.45\% & 100.00 \\
E-n101-k14 & 100 & 14 & 7.1 & 1182.99 & 15 & 0.50\% & 100.00 & 5.28\% & 1183.55 & 14 & 0.38\% & 100.00 \\
P-n45-k5 & 44 & 5 & 8.8 & 537.237$^*$ & 5 &  & 50.56 & 1.98\% & 537.237$^*$ & 5 &  & 4.71$^b$ \\
P-n60-k10 & 59 & 10 & 5.9 & 805.176$^*$ & 11 &  & 7.55 & 3.70\% & 814.436$^*$ & 10 &  & 0.28$^b$ \\
P-n65-k10 & 64 & 10 & 6.4 & 857.481 & 10 & 0.08\% & 100.00 & 1.02\% & 857.481$^*$ & 10 &  & 1.15$^b$ \\
P-n70-k10 & 69 & 10 & 6.9 & 888.060$^*$ & 11 &  & 64.89 & 4.82\% & 912.974 & 10 & 0.48\% & 100.00 \\
X-n101-k25 & 100 & 25 & 4.0 & 30685.8 & 29 & 0.69\% & 100.00 & 9.01\% & 37452.4 & 25 & 3.06\% & 100.00 \\
X-n110-k13 & 109 & 13 & 8.4 & 16734.1 & 14 & 0.18\% & 100.00 & 4.16\% & 17046.5 & 13 & 0.90\% & 100.00 \\
X-n148-k46 & 147 & 46 & 3.2 & 56584.4$^*$ & 50 &  & 13.89 & 1.69\% & 57207.1 & 46 & 0.01\% & 100.00 \\
X-n195-k51 & 194 & 51 & 3.8 & 50040.8 & 58 & 0.56\% & 100.00 & 10.45\% & 58406.0 & 51 & 1.12\% & 100.00 \\
\hline
\end{tabular}}
\captionsetup{width=.88\textwidth}
  \caption*{\footnotesize $^a$ Considering the solution value in BP\&C. $^b$ These instances were solved previously in \S\ref{explit}. The results are replicated.}
  \label{tablelongrun}
\end{table}

Considering both the cases of fixed and unfixed number of routes, 10 additional instances could be solved in the long-run experiments. The relatively small value of stochastic solution obtained in instance $\mathtt{X\text{-}n148\text{-}k46}$ shows that, even in instances with few customers per route, the deterministic equivalent solution can still be a good solution to the stochastic problem.

\section{Conclusions} \label{conclusions}
We considered the VRPSD under optimal restocking, and proposed an algorithm for solving the problem under general distributional assumptions. Computational experiments performed on several literature instances confirmed the effectiveness of the branch-price-and-cut approach. Instances of moderate size (up to 75 nodes) could be solved in reasonable time (up to 5 hours), and larger instances (up to 148 nodes) could be solved in long-runs of the algorithm.

Our algorithm compares favorably with the current approaches in the literature in instances with up to 11 customers per route. In these instances, solving the stochastic problem can be particularly rewarding: the improvements over the deterministic equivalent solution may exceed 10\%. On the other hand, solving the deterministic equivalent problem leads to near-optimal solutions to the stochastic problem when the average number of customers per route is 15 or more, at least under the demand distributions experimented with.

In scenarios where the number of routes is not fixed and the expected demand in a route is not allowed to exceed the vehicle capacity (load factor 1), the optimal solution values under detour-to-depot and optimal restocking are nearly the same. Optimal solutions to the VRPSD often have a number of routes higher than the minimum required for feasibility. This is valid regardless of whether detour-to-depot or optimal restocking is considered. When the average route load approaches the vehicle capacity, opening new routes is a good strategy to reduce expected restocking costs, and in many cases leads to solutions with less transportation (a priori and restocking) costs.

On the other hand, when the number of routes is fixed to the minimum required, the route loads approach the capacity of the vehicle, and the restocking costs become more significant. The additional constraint has an average impact of 2.2\% in the solution value, when applying optimal restocking.

Finally, we performed a range of experiments assuming different load scenarios. The VRPSD has been traditionally studied under the assumption of a fixed load factor of 1, even when optimal restocking is adopted. We showed that by allowing slightly larger load factors (e.g., 1.05), it is possible to improve the objective value by up to 5.2\%. In some cases, it is possible to reduce transportation costs and the number of needed vehicles. Obviously, in practical applications there might be limits on the maximum load factor, otherwise the routes may become too long for complete execution in a single period. Therefore, future research on the VRPSD under optimal restocking should consider load factors larger than 1, and possibly also constraints on the maximum route length or duration.

\bibliographystyle{informs2014} % outcomment this and next line in Case 1
\bibliography{SVRPRefs} % if more than one, comma separated

\end{document}